\documentclass[11pt,a4paper]{amsart}
\usepackage{graphicx,tabularx}
\usepackage{amssymb}
\usepackage{amsfonts}
\usepackage{amsmath}
\usepackage{ifpdf}
\usepackage[all]{xy}
\usepackage{latexsym}
\usepackage{color}

\newtheorem{theorem}{Theorem}[section]
\newtheorem{proposition}[theorem]{Proposition}

\newtheorem{corolario}[theorem]{Corollary}

\newtheorem*{theoremA}{Theorem A}
\newtheorem*{theoremB}{Theorem B}

\theoremstyle{definition}

\newtheorem{example}[theorem]{Example}

\theoremstyle{remark}
\newtheorem{remark}[theorem]{Remark}

\newcommand{\E}{E}

\newcommand{\R}{\mathbb{R}}
\newcommand{\N}{\mathbb{N}}
\newcommand{\cvd}{\ \rule{0.5em}{0.5em}}
\newcommand{\Demo}{{\it Proof: }}

\newcommand{\mR}{g^P}

\newcommand{\Vol}{\mathrm{Vol}}

\newcommand{\opens}{\Omega}
\newcommand{\bopens}{\partial \Omega}
\newcommand{\bifspace}{{\mathcal{M}}}
\newcommand{\biffun}[1]{{\mathcal{F}^{#1}}}
\newcommand{\vol}{\mathcal{V}}

\newcommand{\bord}[1]{\Sigma_{#1}}

\newcommand{\II}{\mathrm I\!\mathrm I}

\begin{document}

\title[Rigidity and bifurcation of CMC hypersurfaces]{Rigidity and bifurcation results for CMC hypersurfaces in warped product spaces}

\author[S. C. Garc\'\i a-Mart\'\i nez]{S. Carolina Garc\'\i a-Mart\'\i nez}
\email{sancarol@ime.usp.br}
\thanks{SCGM was supported by FAPESP (Funda\c{c}\~ao de Amparo \'{a} Pesquisa do Estado de S\~ao Paulo, Brazil) Process 2012/22490-7; MINECO-FEDER project MTM2012-34037 and Fundaci\'on S\'eneca project 04540/GERM/06, Spain. This research is a result of the activity developed within the framework of the Programme in Support of Excellence Groups of the Regi\'{o}n de Murcia, Spain, by Fundaci\'{o}n S\'{e}neca, Regional Agency for Science and Technology (Regional Plan for Science and Technology 2007-2010).
}
\author[J. Herrera]{Jonatan Herrera}
\email{jonatanhf@gmail.com}
\thanks{JH was supported by FAPESP (Funda\c{c}\~ao de Amparo \'{a} Pesquisa do Estado de S\~ao Paulo, Brazil) Process 2012/11950-7; MICINN-FEDER project MTM2010--18099 and grants FQM324, P09-FQM-4496, (J. Andaluc\'ia, the latter also with FEDER funds), Spain.}
\address{Instituto de Matem\'atica e Estat\'istica, Universidade de S\~ao Paulo, Rua do Mat\~ao, 1010 - Cidade Universitaria - S\~ao Paulo, Brazil}
\date{July 28th, 2014}

\begin{abstract}
In this paper, we deduce some rigidity results in warped product spaces under normal variations of CMC hypersurfaces. In particular, we prove the existence of one-parameter families locally rigid on the spatial fiber of Anti-de Sitter Schwarzschild spacetime and one-parameter families with bifurcation points on the spatial fiber of de Sitter Schwarzschild spacetime.
\end{abstract}

\maketitle

\section{Introduction}

Through centuries, the stability of different phenomena has been studied by means of the bifurcation theory. Mathematically, the observed states are modeled as solutions of nonlinear equations, a mathematical context where the intuitive notion of stability is made precise. Roughly speaking, one expects that a slight change of a parameter in a system should not have great influence, but rather that stable solutions change continuously in a unique way. That expectation is verified by the Implicit Function Theorem. Consequently, as long as a continuous branch of solutions preserves its stability, no dramatic change is observed when the parameter is varied. However, if that \emph{ground state} loses its stability when the parameter reaches a critical value, then the state is no longer observed, and the system itself organizes a new stable state that \emph{bifurcates} from the ground state.

In \cite{S-W}, Smoller and Wasserman present a simple criterion based on the Index Morse determining the existence of bifurcation points which, together with the Implicit Function Theorem, has enabled several authors to study different stability problems. For instance, Bettiol and Piccione \cite{BP} used this technique to obtain results about the existence and non-existence of metrics with constant scalar curvature and as close as desired to a given homogeneous metric on the sphere. Jin, Lin and Xu in \cite{JLX} obtained some multiplicity results for the Yamabe equation on the sphere $\mathbb{S}^n$. Recently, Ram\'\i rez-Ospina \cite{RO} proved the existence of uncountably many unit volume constant scalar curvature metrics in product manifolds $\mathbb{T}^k\times M$ where $\mathbb{T}^k$ is a flat $k$-torus and $M$ is a compact manifold. Among all these works, we point out one authored by Al\'\i as and Piccione \cite{AP}, where it is shown the existence of an infinite sequence of isometric embeddings of tori with constant mean curvature (CMC) in Euclidean spheres that are not isometrically congruent to the CMC Clifford torus.

Motivated by this work, we propose to study the rigidity of families of CMC hypersurfaces in different (and more general) contexts. Roughly speaking, a family $\mathcal F$ of CMC hypersurfaces is (locally) rigid if given any other CMC hypersurface $\Sigma$ which is sufficiently close to some element of $\mathcal F$, then $\Sigma$ must be \emph{congruent} to some element of $\mathcal F$ (see Section~\ref{sub:niftheory} for details). Constant mean curvature hypersurfaces have a well-known variational characterization as critical points of the area functional under variations preserving certain volume. Our aim is to show when a family of hypersurfaces is rigid or it has some bifurcation point. A natural occurrence of foliations by CMC hypersurfaces appears in certain warped product spaces, which are the main object of our study.

Finally, we will apply our results to the spatial fiber of the so-called de Sitter and Anti-de Sitter Schwarzschild spacetimes. These classical solutions of the vacuum Einstein equation for static, spherically symmetric spacetimes describe the simplest model of a black hole which is known. Recently, these spaces have received an increasing attention, especially with the so-called AdS/CFT correspondence (or Maldacena duality), a conjectured equivalence between a string theory on a {\em bulk} space (typically, the product of Anti-de Sitter spaces by round spheres) and a quantum field theory without gravity on the {\em boundary} of the initial space, which behaves as an hologram of lower dimension (see \cite{fisicos} for a detailed introduction on the topic).

The paper is organized as follows. In Section 2, we establish the basic notions of bifurcation and local rigidity; and we present the variational setting in which we will work.
The bifurcation phenomenon under volume-preserving variations of CMC hypersurfaces is analyzed in Section 3. As a first step, we consider the case of families of CMC hypersurfaces whose associated potential function of the Jacobi operator is constant (see \eqref{potential}). For this general setting, including the space of CMC hypersurfaces which are orbits for some isometric action of a Lie group (as it happens in the warped product case), we obtain simple sufficient conditions to guarantee the stability of an one-parameter family of hypersurfaces. Our main results are summarized in the following statements (see Theorem \ref{teo1}, Corollary \ref{cor1} and Remark \ref{hypopen} (i) for details).

\begin{theoremA}
Let $\{\partial\Omega_\gamma\}_{\gamma\in I}$ be a family of CMC hypersurfaces in an $n$-dimensional Riemannian manifold $(M,g)$. Assume that, for each $\gamma$ the map  
\[Q(\gamma):=||\II_{\gamma}||^2+(n-1)\mathrm{Ric}(n_\gamma,n_\gamma)  \]
is constant on $\partial \Omega_\gamma$, where $||\II_{\gamma}||^2$ is the squared norm of the second fundamental of $\bopens_\gamma$
and $\mathrm{Ric}(n_\gamma,n_\gamma)$ is the normalized Ricci curvature of $M$ evaluated at $n_\gamma$ the unit normal to $\partial\Omega_\gamma$.

Then, the family $\{\bopens_{\gamma}\}_{\gamma\in I}$ is rigid for all $\gamma$ if
\[\mu_1(\gamma)-Q(\gamma)>0\]
where $\mu_1(\gamma)$ represents the first nonzero eigenvalue of the Laplacian defined on $\bopens_\gamma$ with the induced metric. In particular, such a family is rigid if one of the following statements hold:

\begin{itemize}
\item[$i)$] $\bopens_\gamma$ is convex\footnote{We recall that $\bopens_\gamma$ is convex if and only if its second fundamental form is positive definite with respect to the inward normal direction.} and has non-positive Ricci curvature,
\item[$ii)$] $\mathrm{Ric}(n_\gamma,n_\gamma)\leq-1/(n-1)||\II_\gamma||^2$,
\item[$iii)$] $\bopens_\gamma$ is Ricci flat and $\mu_1(\gamma)\geq||\II_\gamma||^2$, where $\mu_1(\gamma)$ is the first nonzero eigenvalue of $\Delta$ on $\bopens_\gamma$.
\end{itemize}

\end{theoremA}

\smallskip

For the particular case of warped product spaces, we give sufficient conditions for both the existence and non-existence of bifurcation points. The results for this case are summarized as follow (see Propositions \ref{prop2.3}, \ref{prop5.1} for more details).

\begin{theoremB}
Let $((r_1,r_2)\times P,dr^2+\alpha^2(r)g^P)$ be a warped product space where $P$ is an $(n-1)$-dimensional closed (compact and without boundary) manifold; and consider the one-parameter family of CMC hypersurfaces $\{\{r\}\times P\}_{r\in(r_1,r_2)}$.

\begin{itemize}
\item[(i)] If
\[\dot{\alpha}^2(r)-\ddot{\alpha}(r)\alpha(r)<\frac{\hat{\mu}_1}{n-1}\] for all $r\in (r_1,r_2)$ where $\hat{\mu}_1$ is the first non-zero eigenvalue of the Laplacian on $(P,g^P)$, then $\{\{r\}\times P\}_{r\in(r_1,r_2)}$ has no bifurcation point.\smallskip

\item[(ii)] If
\[\lim_{r\rightarrow r_2 }\dot{\alpha}^2(r)-\ddot{\alpha}(r)\alpha(r)=\infty
\]
then, there exists an infinite sequence of bifurcation points.
\end{itemize}
\end{theoremB}
Furthermore, we give some geometric conditions (on terms of the scalar and mean curvatures of warped product and
the scalar curvature of its fiber), which allow us to ensure the rigidity of an one-parameter family of 
CMC hypersurfaces in a warped product, see Corollary \ref{corsc}. This type of situation occurs in General Relativity, and we will discuss in Section 4 an application of these results to the spatial fiber of the de Sitter and Anti-de Sitter Schwarszchild spacetimes. We show that the family of spheres around the blackhole are locally rigid for the former (including the classical Schwarszchild spacetime) while the collapsing spheres on the latter have infinitely many bifurcation points.

\section{Preliminaries}\label{section2}

Let us state some of the basic elements, results and notations that we are going to use in the rest of the paper. Along this paper, we will consider an $n$-dimensional Riemannian manifold $(M^n,g)$. In general, a local coordinate system will be denoted by $(x_{i})_{i=1}^{n}$ and, in the adapted case for the boundary, we will assume that $\partial/\partial x_n$ is the normal vector to $\partial M$ pointing inward. The volume and area elements of $M$ and $\partial M$ will be denoted by $dv$ and $d\sigma$ respectively. \\

\subsection{On bifurcation theory}
\label{sub:niftheory}
This section is devoted to introduce roughly the basic framework about bifurcation theory needed for the rest of this paper. For details, we refer the readers to \cite{S-W,AP} and the references therein.

Let us consider a Banach space $\bifspace$ and a continuous path of $C^{k}$-functionals (with $k\geq 2$) $\biffun{\lambda}:\bifspace\rightarrow \R$ where $\lambda$ varies on a prescribed interval $I$. Assume that we also have a continuous path of critical points $\{x_{\lambda}\}\subset \bifspace$ of the corresponding functionals, that is, $d\biffun{\lambda}(x_{\lambda})=0$.

Fix $\lambda_*\in I$. We will say that $\lambda_*$ is \textit{a point of bifurcation} if there exists a sequence $\{\lambda_n\}_n\subset I$ and a sequence $\{x_n\}_n$ satisfying:
\begin{itemize}
\item[(1)] $\lim\limits_{n\rightarrow \infty}\lambda_n=\lambda_*$ and $\lim\limits_{n\rightarrow \infty}x_n=x_*$,
\item[(2)] $d\biffun{\lambda_n}(x_n)=0$ for all $n$,
\item[(3)] $x_n\neq x_{\lambda_n}$ for all $n$.
\end{itemize}
If $\lambda_*$ is not a bifurcation point, we will just say that the family $\{x_{\lambda}\}_{\lambda}$ is \textit{locally rigid} at $\lambda_*$.

One of the classical criterion to determine when a point is of bifurcation is related with the so-called \textit{Morse index}. Recall that the Morse index of a critical point $x_{0}$  of a smooth functional $\biffun{}$ is equal to the dimension of the maximal subspace of the tangent space $T_{x_{0}}\bifspace$ where the second variation $d^2(\biffun{})_{x_{0}}$ is negative definite. Such an index will be denoted in general by $i(x_{0},\biffun{})$.

On the one hand, essentially, a variation of the Morse index $i(\biffun{\lambda},x_{\lambda})$ along the interval $I$ will indicate the existence of a bifurcation point (see for instance \cite{S-W} for details). More precisely, under suitable fredholmness assumptions, if there exist $\lambda_1,\lambda_2\in I$ with $\lambda_1<\lambda_2$ such that: $d^2(\biffun{\lambda_i})_{x_{\lambda_i}}$ are non-singular for $i=1,2$ and $i(x_{\lambda_1},\biffun{\lambda_1})\neq i(x_{\lambda_2},\biffun{\lambda_2})$, then there exists a bifurcation point $\lambda_{*}\in (\lambda_1,\lambda_2)$.

On the other hand, using the Implicit Function Theorem, one sees that if $d^2(\biffun{\lambda_*})_{x_{\lambda_*}}$ is non singular then the family $\{x_{\lambda}\}$ is \textit{locally rigid} at $\lambda_*$. In particular, if $i(\biffun{\lambda},x_{\lambda})=0$ for all $\lambda$ in $I$, then there is no bifurcation point.

In this paper, using the above conditions, we will study local rigidity and bifurcation by analyzing the spectrum of $d^2(\biffun{\lambda})_{x_{\lambda}}$ for all $\lambda$. Essentially, we will determine the number of negative eigenvalues for each $\lambda$ (counting its multiplicity) and we will study the evolution of such a number as $\lambda$ runs along $I$.

\subsection{Stating the variational problem}

Let us recall the classical variational setup of the constant mean curvature (CMC) hypersurfaces as critical points of the area functional under variations preserving the volume. 

\smallskip

In order to use the bifurcation theory above mentioned, let us describe the elements involved in this approach. The Banach space $\bifspace$ will be the space of open subsets $\opens$ of an $n$-dimensional Riemannian manifold $(M,g)$ with compact closure and whose smooth boundary is a connected and orientable smooth manifold. The variations on $\bifspace$ will be defined in the following way: a differentiable map $X:(-\epsilon,\epsilon)\times \bopens \rightarrow M$ is called a {\em CMC hypersurface variation of $\bopens$} if $(i)$ $X_{t}:\bopens\rightarrow M$, defined by $X_{t}(x)=X(t,x)$ is an immersion whenever $|t|<\epsilon$, for $\epsilon>0$ small enough; and $(ii)$  $X(0,x)=i(x)$ where $i$ is the inclusion map.

Given an element $\opens\in \bifspace$ and a variation $X$ of $\bopens$ denote by $\bopens_{t}=X_t(\bopens)$. For values of $t$ small enough, $\bopens_{t}$ is also a connected and oriented smooth submanifold. Moreover, it contains an open subset $\opens_{t}$ whose closure is also compact. In conclusion, the hypersurface variation defines in a natural way a variation of the open subset $\opens$ denoted by $\opens_{t}$, which is also an element of $\bifspace$. Among all the hypersurface variations, we are going to be interested in the so-called \textit{normal} variations. A hypersurface variation $X$ is normal if the variation vector field $\frac{\partial X}{\partial t}\big|_{t=0}$ is parallel to the unit normal vector field $N$ with respect to $\Omega$ and the volume is preserved. In particular, if we denote by $T_{\opens}^N \bifspace$ the linear subspace of $T_{\opens} \bifspace$ only determined by normal variations, we can deduce that $T_{\opens}^N \bifspace$ is naturally identified with the set of smooth functions with compact support $\mathcal{C}^\infty_0(\bopens)$, where such an identification is essentially obtained by defining  $f=\left\langle \frac{\partial X}{\partial t}\big|_{t=0}, N\right\rangle\in\mathcal{C}_0^\infty(M)$ for a normal variation $X$. For the volume-preserving property, let us recall that the volume form varies as

\begin{equation}\label{forvol}
\Vol(\opens_{t})= \int_{\opens_{t}}dv=\int_{\opens}dv+ \int_{[0,t]\times \bopens}X^*(dv).
\end{equation}

\noindent and so, the normal variations has an associated function $f$ satisfying

\begin{align}\label{volpre}\int_{\bopens} f d\sigma=0.
\end{align}

Finally, we will consider the area functional restricted to such a volume-preserving variations. The Lagrange multiplier method leads us then to the functional $\biffun{\lambda}$ defined as:

\begin{equation}
\begin{array}{cccc}
\biffun{\lambda}: & \bifspace & \longrightarrow & \R\\
 & \opens & \longmapsto & \textrm{Area}(\bopens) + \lambda \Vol(\opens),\end{array}
\end{equation}
where $\lambda$ varies on a prescribed interval $I$ (to be determined later).

\smallskip

Now, we are ready to determine the first variation of our functional $\biffun{\lambda}$. For this, consider $X$ a normal variation and $f$ its associated function satisfying \eqref{volpre}. Denote by

\[
A(t)=\textrm{Area}(\bopens_t),\quad \vol(t)=\Vol(\opens_{t}).
\]

Taking into account \eqref{forvol}, we deduce that

\[
\frac{\partial}{\partial t}\Big|_{t=0}A(t)=-(n-1)\int_{\bopens}fH d\sigma\quad \hbox{\and} \quad \frac{\partial}{\partial t}\Big|_{t=0} \vol(t)=\int_{\bopens} f d\sigma
\] where $H$ is the mean curvature on $\partial \Omega$ (for more details see \cite[Lemma 2.1]{BCE}).

So, that the first variation of our functional takes the following form

\begin{align}
\delta\biffun{\lambda}_{\opens}(f)=& \frac{\partial}{\partial t}\Big|_{t=0} A(t) + \lambda \frac{\partial}{\partial t}\Big|_{t=0}\vol(t)\nonumber\\
= &\int_{\bopens} \left(-(n-1)H + \lambda\right)f d\sigma. 
\end{align}

Therefore, the critical points of $\biffun{\lambda}$ under normal variations are open subsets of $M$ such that its boundaries are CMC-hypersurfaces with $H=\lambda/(n-1)$.

\smallskip

\begin{remark}\label{hypopen}\label{hypboundos}(i) It is important to remark that, for these computations, there is no real dependence on the open set $\Omega$ but on the hypersurface $\partial\Omega$. In fact, in the literature, it is more common to work in terms of hypersurfaces. However, for simplicity on this case, we have preferred this approach of open sets.

(ii) Note also that this is also valid for a more general setting. Assume that $\bifspace$ is the space of open subsets $\opens\subset M$ whose boundary is union of two disjoint sets $\partial\Omega=\Sigma_1\cup\Sigma_2$. We will assume that one of them, $\Sigma_1$, is a fixed set and so that the hypersurface variations only affects $\Sigma_2$. Under this assumption, the critical points of the functional will be open subsets $\Omega$ such that their boundaries are union of a (fixed) set $\bord{1}$ and a CMC hypersurface $\bord{2}$.
\end{remark}

As we are interested in rigidity results, we need to study the second variation of the functional on the critical points or, more precisely, its spectrum. Assume that $\opens$ is a critical point for the functional, that is, that the boundary $\bopens$ is an $(n-1)$-dimensional hypersurface with constant mean curvature $H=\lambda/(n-1)$. The second variation of the functional on $\opens$ has been already computed in \cite{BCE}, and it has the following form:
\begin{equation}\label{secondvariationsur}
 \delta^2 \biffun{\lambda}_{\opens}(f):=\Delta f - \left(||\II||^2+(n-1)\mathrm{Ric}(n_0,n_0)\right)f
 \end{equation}

\noindent where $\Delta$ represent the non-negative Laplace-Beltrami operator and $||\II||^2$ the squared norm of the second fundamental form both in $\bopens$ with the induced metric by $g$; $\mathrm{Ric}(n_0,n_0)$ is the normalized Ricci curvature of $M$ evaluated at $n_0$ the normal vector to $\opens$ pointing outward and $f$ satisfies \eqref{volpre}.

\section{Bifurcation under cmc hypersurface variations}
\subsection{General bifurcation results}
Consider now an one-parameter family $\{\opens_\gamma\}_\gamma$ of open subsets in $M$ such that the boundary of each $\opens_\gamma$, denoted by $\bopens_{\gamma}$, is a compact hypersurface with constant mean curvature $H(\gamma)$. In particular, each element $\opens_{\gamma}$ is a critical point for the functional $\biffun{H(\gamma)(n-1)}$ where we can study its rigidity character. As our aim is to see if there is a change in the Morse index, we need to
study the eigenvalues of the linear map $\delta^2(\biffun{H(\gamma)(n-1)})_{\opens_\gamma}(f)$.
At this point, in order to avoid confusion, let us introduce the following notation 
\[
\lambda(\gamma):=H(\gamma)(n-1),\qquad {\mathcal A}_{\gamma}:=\delta^2(\biffun{H(\gamma)(n-1)})_{\opens_\gamma}.
\]
 Notice, if the potential function in \eqref{secondvariationsur} is constant then the eigenfunctions of such a map will coincide with the eigenfunctions of the Laplacian and, even more, if we denote such potential by 
\begin{eqnarray}\label{potential}
Q(\gamma):=||\II_\gamma||^2+(n-1)\mathrm{Ric}(n_{\gamma},n_{\gamma})
\end{eqnarray}
 and if $f$ is an eigenfunction of the Laplacian associated to an eigenvalue $\mu(\gamma)$, then $\overline{\mu}(\gamma)$ is an eigenvalue for the second variation of the functional, where

\[
\overline{\mu}(\gamma)=\mu(\gamma)-Q(\gamma)
\]
(here $\gamma$ represents its dependence respect $\opens_{\gamma}$). Finally, by the spectral theorem \cite[Theorem A.I.4]{BGM} we know that all the eigenvalues of the (non-negative) Laplacian $\Delta$ on $\bopens_{\gamma}$, are determined by a sequence of eigenvalues $\{\mu_i(\gamma)\}_{i\in \mathbb{N}_0}$ satisfying\footnote{The set $\mathbb{N}_0=\N\cup\left\{0\right\}=\left\{0,1,2,\cdots\right\}$.}:
\[\mu_0(\gamma)=0<\mu_1(\gamma)\leq\dots\leq\mu_{i}(\gamma)\leq \mu_{i+1}(\gamma)\leq\dots\quad \textrm{and}\quad \lim_{i}\mu_i(\gamma)=\infty,\] repeated according to their multiplicity.

\begin{remark}
Since the potential $Q(\gamma)$ is constant, if $f$ is an eigenfunction for $\Delta$ also it is for the second variation of the functional ${\mathcal A}_\gamma$ (see \eqref{secondvariationsur}). However, for our problem, an eigenfunction for the latter should additionally satisfy the integral condition \eqref{volpre}. In particular, the eigenfunctions associated to the eigenvalue zero of $\Delta$, which are non-zero constants, are not valid eigenfunctions for $\mathcal{A}_{\gamma}$, as they do not satisfy such an integral condition. In conclusion, the eigenvalues on (\ref{valfun}) have to be considered for $i\in\mathbb{N}$.
\end{remark}

In conclusion, assuming that the potential is constant, all the eigenvalues of ${\mathcal A}_{\gamma}$ have the following form:
\begin{align}\label{valfun}
\overline{\mu_i}(\gamma)=\mu_i(\gamma)-Q(\gamma),\quad\textrm{for every }i\in\mathbb{N}.
  \end{align}
As we have determined the spectrum of the second derivative completely, we can characterize both, the rigidity and the existence of bifurcation points on the family $\{\opens_{\gamma}\}_{\gamma}$. Let us start giving a simple sufficient condition for the rigidity:

\begin{theorem}\label{teo1}
Let $\{\opens_{\gamma}\}_{\gamma\in I}$ be a family of open subsets on $M$ whose boundary $\bopens_{\gamma}$ satisfies both: it is a hypersurface with constant mean curvature $H(\gamma)$ and its associated potential $Q(\gamma)$ is constant on $\bopens_{\gamma}$. Then, such a family is rigid if for all $\gamma$,
\[\mu_1(\gamma)-Q(\gamma)>0\]
where $\mu_1(\gamma)$ represents the first nonzero eigenvalue of the Laplacian defined on $\bopens_\gamma$.

\end{theorem}

\Demo As we have already observed, $\opens_\gamma$ is a critical point for the functional $\biffun{\lambda(\gamma)}$. Moreover, we have also proved that all the eigenvalues of ${\mathcal A}_{\gamma}$ are determined by \eqref{valfun}. So, taking into account that, for all $i\in\mathbb{N}$

\[
\overline{\mu_i}(\gamma)=\mu_i(\gamma)-Q(\gamma)\geq \mu_1(\gamma)-Q(\gamma)>0,
\]

then the second derivative of the functional is non singular for all $\gamma$. Therefore, we obtain the desired result.
\cvd

\smallskip

In particular, we obtain the following simple conditions to ensure rigidity:
\begin{corolario}\label{cor1}
Consider a family $\{\opens_{\gamma}\}_{\gamma}$ of open subsets of $M$ whose boundary $\partial\Omega_\gamma$ is a CMC hypersurface. Such a family is rigid if one of the following conditions holds:
\begin{itemize}
\item[$i)$] $\bopens_\gamma$ is convex and has non-positive Ricci curvature,
\item[$ii)$] $\mathrm{Ric}(n_0,n_0)\leq-1/(n-1)||\II_\gamma||^2$,
\item[$iii)$] $\bopens_\gamma$ is Ricci flat and $\mu_1(\gamma)\geq||\II_\gamma||^2$, where $\mu_1(\gamma)$ is the first nonzero eigenvalue of $\Delta$ on $\bopens_\gamma$.
\end{itemize}
\end{corolario}

\smallskip

Now, let us formalize the criterion that we will use to ensure the existence of bifurcation points (recall Section \ref{section2}).

\begin{theorem}\label{teo2}
Let $\{\opens_{\gamma}\}_{\gamma\in I}$ be a family of open subsets on $M$ whose boundary $\bopens_{\gamma}$ is a hypersurface with constant mean curvature $H(\gamma)$. If there exist two values $\gamma_1,\gamma_2\in I$ with $\gamma_1<\gamma_2$ satisfying: 
\begin{itemize}
\item[(i)] $\overline{\mu_{i}}(\gamma_j)\neq 0$ for all $i$ and $j=1,2$;
\item[(ii)] there exists $i_0$ such that $\overline{\mu_{i_0}}(\gamma_1)\overline{\mu_{i_0}}(\gamma_2)<0$,
\end{itemize} then there exists a bifurcation point $\gamma_*\in (\gamma_1,\gamma_2)$.
\end{theorem}

\Demo 
For this result we have to show both, that ${\mathcal A}_{\gamma}$ is non-singular for $j=1,2$; and that $i(\opens_{\gamma_1},\biffun{\lambda(\gamma_1)})\neq i(\opens_{\gamma_2},\biffun{\lambda(\gamma_2)})$. The former condition is a direct consequence of assumption (i), as all the eigenvalues of such a linear map are nonzero. For the latter, observe that hypothesis (ii) is ensuring that the eigenvalue $\overline{\mu_{i_0}}$ is changing the sign between $\gamma_1$ and $\gamma_2$. Moreover, as the eigenvalues $\overline{\mu_i}$ are ordered, we can ensure that the number of negative eigenvalues between $\gamma_1$ and $\gamma_2$ has changed, so $i(\opens_{\gamma_1},\biffun{\lambda(\gamma_1)})\neq i(\opens_{\gamma_2},\biffun{\lambda(\gamma_2)})$.
Then the result is obtained. \cvd
\\
\smallskip

\subsection{CMC hypersurfaces bifurcation in warped product spaces}
\label{hybifwsection}
Let us assume now that $(M,g)$ is an $n$-dimensional warped product space, that is, $M=(r_1,r_2)\times P$ is a product manifold endowed with the following metric
\begin{equation}\label{eq1}
g=dr^2 +\alpha^2(r)\mR
\end{equation}
where $(P,\mR)$ is an $(n-1)$-dimensional closed Riemannian manifold (i.e. a compact manifold without boundary) and $\alpha$ is a smooth positive function on $(r_1,r_2)$.

 These spaces naturally define a family of open subsets which can be realized as critical points of the area-volume functional for some special $\lambda$. Such a family, denoted by $\{\opens_{r}\}_{r\in (r_1,r_2]}$, is formed by elements $\opens_{r}=(r_1,r)\times P$ whose boundary of each $\opens_r$ is composed by a fixed set $\bord{1}=\{r_1\}\times P$ and other set $\bord{2}=\{r\}\times P$. 
 
So, recalling Remark \ref{hypboundos} (ii) we have that $\opens_r$ is a critical point for the functional $\biffun{\lambda(r)}$ where $\Sigma_2$ is a hypersurface with constant mean curvature (with respect to the inward unit normal $-\partial r$) given by\footnote{By notation, the dot will denote derivative of a real function.} \[H(r)=-(n-1)\frac{\dot{\alpha}(r)}{\alpha(r)}.\]

 Now our aim is to study at what extent such a family is rigid. The particular advantage of the warped product spaces is that we can explicitly determine both, $||\II||^2$ and $\mathrm{Ric}_{p_0}$ in terms of the warping function $\alpha$. In fact, by simple computations we obtain that

 \[
 ||\II||^2=(n-1)\left(\frac{\dot{\alpha}}{\alpha}\right)^2\quad\textrm{and} \quad \mathrm{Ric}_{p_0}(\partial r, \partial r)=-\frac{\ddot{\alpha}\alpha}{\alpha^2}
 \]
 which shows that the potential $Q(r)$ is constant on the boundary of the open set and leads us to the following expression for the second variation (recall (\ref{secondvariationsur}) and \eqref{potential})

\begin{equation}\label{eqbif}\begin{array}{rl}
{\mathcal A}_{r}(f)= & \Delta f - (n-1)\left(\dfrac{\dot{\alpha}^2-\ddot{\alpha}\alpha}{\alpha^2}\right)f.
\end{array}\end{equation}

\smallskip

Moreover, the induced metric on $\bord{2}$ is just $g^P$ multiplied by the constant $\alpha^2(r)$, which allows us to relate easily the spectrum of $\Delta$ with the corresponding spectrum of the Laplacian on $(P,g^P)$. In fact, if $\hat{\mu}$ is an eigenvalue of the latter, $\mu=\hat{\mu}/\alpha^2(r)$ will be an eigenvalue for the former.

Therefore, as $P$ is a compact manifold, applying again the spectral theorem we deduce that the eigenvalues of the Laplacian on $(P,g^P)$ are determined by a sequence $\{\hat{\mu}_i\}_{i\in\mathbb{N}_0}$ satisfying:

\[\hat{\mu}_0=0<\hat{\mu}_1\leq\dots\leq\hat{\mu}_{i}\leq \hat{\mu}_{i+1}\leq\dots\quad \textrm{and}\quad \lim_{i}\hat{\mu}_i=\infty,\] repeated according to their
multiplicity.

Summarizing, we deduce that the eigenvalues of the linear map ${\mathcal A}_{r}$ have the  form described on (\ref{valfun}) with

\begin{equation}\label{valfun1}\mu_i(r)=\frac{\hat{\mu}_i}{\alpha^2(r)}\quad \hbox{and}\quad Q(r)=\frac{(n-1)(\dot{\alpha}^2(r)-\ddot{\alpha}(r)\alpha(r))}{\alpha^2(r)}
\end{equation}

\noindent where $\hat{\mu}_i$ are the eigenvalues of the Laplacian of $(P,g^P)$. \\

Then, we can just give simple conditions for both, rigidity and existence of bifurcation points in terms of the warping function $\alpha$ and the eigenvalues of the Laplacian on $(P,g^P)$. For rigidity, we will just give a simple translation of Theorem \ref{teo1} (the proof is trivial from such a theorem and \eqref{valfun1}):

\begin{proposition}\label{prop2.3}
Let $\left((r_1,r_2)\times P,dr^2+\alpha^2(r)\mR \right)$ be a warped product space where $P$ is an $(n-1)$-dimensional closed manifold and consider the one-parameter family $\{\opens_r\}_r$ defined as above. If \[\dot{\alpha}^2(r)-\ddot{\alpha}(r)\alpha(r)<\frac{\hat{\mu}_1}{n-1}\] for all $r\in (r_1,r_2)$, then $\{\opens_r\}_r$ has no bifurcation points.
\end{proposition}

\begin{example}
From here, we easily conclude that the pseudo hyperbolic space (for more details about this space, see \cite{T}) \[\left(
(r_1,r_2)\times\mathbb{S}^{n-1}, dr^2+\,e^{2r}\,g^{\mathbb{S}^{n-1}}\!\right)  \]
is foliated by a locally rigid family of CMC hypersurfaces (recall Remark \ref{hypopen} (i)). The slices $\{r\}\times\mathbb{S}^{n-1}$ determine an one-parameter family of CMC hypersurfaces without bifurcation points 
because $\dot{\alpha}^2(r)-\ddot{\alpha}(r)\alpha(r)=0$ and $\hat{\mu}_1=n-1$ (see \cite[Proposition C.I.1]{BGM}).\\
\end{example}

\smallskip

Next, our objective is to give conditions for $\alpha$ which ensure the existence of bifurcation points of the family $\{\opens_{r}\}_{r\in(r_1,r_2)}$. Fix a value $r_0\in (r_1,r_2)$ and observe that, even if we do not know the sign on $\overline{\mu_i}(r_0)$, we can ensure that there exists $i_0\in \mathbb{N}$ such that for all $i\geq i_0$, $\overline{\mu_i}(r_0)>0$ (recall that $\hat{\mu}_i\rightarrow \infty$, \eqref{valfun} and \eqref{valfun1}). So, we only need to impose conditions which ensure that, for some $r>r_0$, $\overline{\mu_{i_0}}(r)<0$. This idea is exploited to obtain the following result.
\begin{proposition}\label{prop5.1}
Let $\left((r_1,r_2)\times P,dr^2+\alpha^2(r)\mR  \right)$ 	be a warped product space where $P$ is an $(n-1)$-dimensional closed manifold and consider the one-parameter family of open subsets $\opens_{r}=(r_1,r)\times P$ whose boundary is composed by a fix set $\bord{1}=\{r_1\}\times P$ and a CMC hypersurface $\bord{2}=\{r\}\times P$ (i.e., critical points for the functional $\biffun{\lambda(r)}$).
If
\begin{equation}\label{cond1}
\lim_{r\rightarrow r_2}\dot{\alpha}^2(r)-\ddot{\alpha}(r)\alpha(r)=\infty
\end{equation}
then, there exist infinitely many bifurcation points.
\end{proposition}
\Demo 
Let us denote by \[h(r):=(n-1)\left(\dot{\alpha}^2(r)-\ddot{\alpha}(r)\alpha(r)\right)\] so that we can write the eigenvalues $\overline{\mu_{i}}(r)$ in the following way (recall \eqref{valfun} and \eqref{valfun1}):

\[
\overline{\mu_i}(r)=\frac{1}{\alpha^2(r)}\left(\hat{\mu}_i-h(r)\right).
\]

Fix $r_0\in (r_1,r_2]$ and consider $i_0\in \N$ such that $\overline{\mu_{i_0}}(r_0)>0$, that is, with
$ \hat{\mu}_{i_0}>h(r_0)$
(which is possible because $\lim_{i}\hat{\mu}_{i}=\infty$). As $h$ tends to infinity when $r\rightarrow r_2$ (recall (\ref{cond1})), we can find some value $\widetilde{r}>r_0$ such that
$\hat{\mu}_{i_0}<h(\widetilde{r})$,
and so, with $\overline{\mu_{i_0}}(\widetilde{r})<0$. Now, as $h$ is a continuous function and satisfies

\[h(r_0)<\hat{\mu}_{i_0}<h(\widetilde{r})
\]
we can obtain two values $r'_0,r'_1$ with $r_0\leq r'_0<r'_1\leq \widetilde{r}$ such that 

\[
\hat{\mu}_{i_{-1}}<h(r'_0)<\hat{\mu}_{i_0}<h(r'_1)<\hat{\mu}_{i_1}
\] where $\hat{\mu}_{i_{-1}}$ (resp. $\hat{\mu}_{i_{1}}$) denotes the largest (resp., smallest) eigenvalue of the Laplacian $\Delta$ less (resp., greater) than $\hat{\mu}_i$. Then we deduce both, that $\overline{\mu_{i_0}}(r'_0)\overline{\mu_{i_0}}(r'_1)<0$ and that $\overline{\mu_{i}}(r'_j)\neq 0$ for all $i\in\mathbb{N}_0$ and $j=1,2$. Therefore, Theorem \ref{teo2} ensures the existence of a bifurcation point in the interval $(r_0,\widetilde{r})$. Now, we can repeat the process with $\widetilde{r}$, obtaining another bifurcation point in the interval $(\widetilde{r},r_2)$ for some $r_2>\widetilde{r}$. In conclusion, by induction, we prove the existence of infinitely many points of bifurcation.
\cvd

\begin{figure}
\centering
\ifpdf
  \setlength{\unitlength}{1bp}%
  \begin{picture}(235.19, 189.06)(0,0)
  \put(0,0){\includegraphics{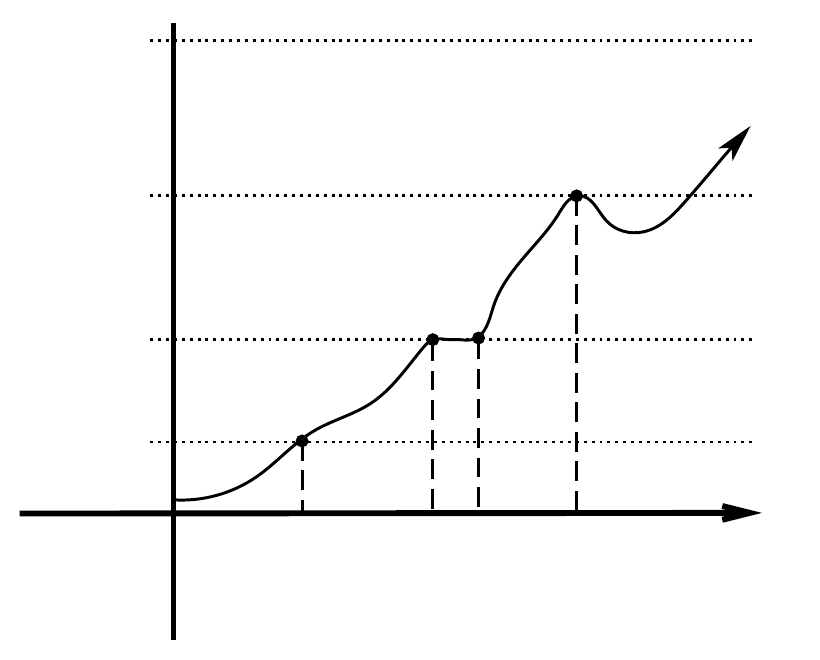}}
  \put(30.88,57.36){\fontsize{8.54}{10.24}\selectfont $\hat{\mu}_{i}$}
  \put(24.38,87.88){\fontsize{8.54}{10.24}\selectfont $\hat{\mu}_{i+1}$}
  \put(24.38,130.22){\fontsize{8.54}{10.24}\selectfont $\hat{\mu}_{i+2}$}
  \put(24.38,175.02){\fontsize{8.54}{10.24}\selectfont $\hat{\mu}_{i+3}$}
  \put(205.51,33.73){\fontsize{8.54}{10.24}\selectfont $r$}
  \put(51.95,34.22){\fontsize{8.54}{10.24}\selectfont $r_1$}
  \put(212.51,153.73){\fontsize{11.38}{13.66}\selectfont $h$}
  \put(81.82,34.76){\fontsize{5.69}{6.83}\selectfont $r_{A}$}
  \put(118.06,34.76){\fontsize{5.69}{6.83}\selectfont $r_{B}$}
  \put(133.59,34.76){\fontsize{5.69}{6.83}\selectfont $r_{C}$}
  \put(162.30,34.76){\fontsize{5.69}{6.83}\selectfont $r_{D}$}
  \end{picture}%
\else
  \setlength{\unitlength}{1bp}%
  \begin{picture}(235.19, 189.06)(0,0)
  \put(0,0){\includegraphics{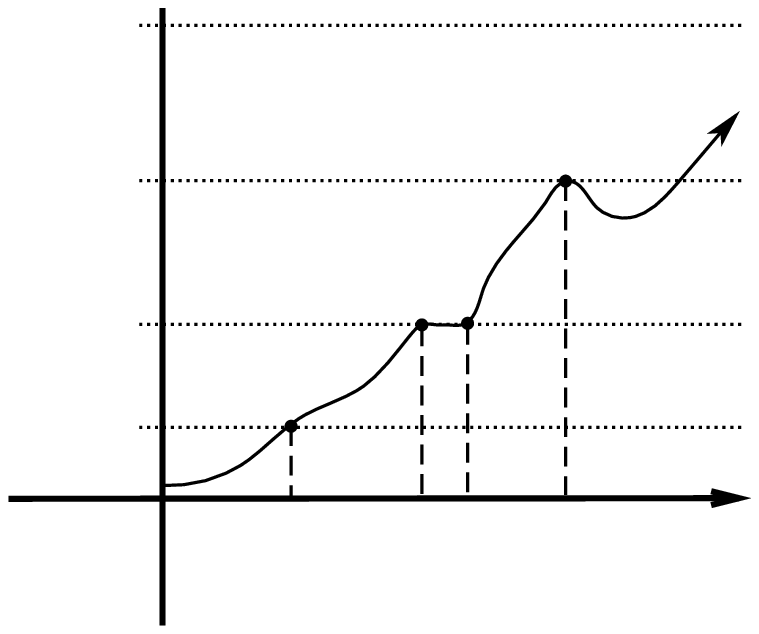}}
  \put(30.88,57.36){\fontsize{8.54}{10.24}\selectfont $\hat{\mu}_{i}$}
  \put(24.38,87.88){\fontsize{8.54}{10.24}\selectfont $\hat{\mu}_{i+1}$}
  \put(24.38,130.22){\fontsize{8.54}{10.24}\selectfont $\hat{\mu}_{i+2}$}
  \put(24.38,175.02){\fontsize{8.54}{10.24}\selectfont $\hat{\mu}_{i+3}$}
  \put(205.51,33.73){\fontsize{8.54}{10.24}\selectfont $r$}
  \put(51.95,34.22){\fontsize{8.54}{10.24}\selectfont $r_1$}
  \put(212.51,153.73){\fontsize{11.38}{13.66}\selectfont $h$}
  \put(81.82,34.76){\fontsize{5.69}{6.83}\selectfont $r_{A}$}
  \put(118.06,34.76){\fontsize{5.69}{6.83}\selectfont $r_{B}$}
  \put(133.59,34.76){\fontsize{5.69}{6.83}\selectfont $r_{C}$}
  \put(162.30,34.76){\fontsize{5.69}{6.83}\selectfont $r_{D}$}
  \end{picture}%
\fi
\caption{\label{fig1}The function $h$ let us deduce where we can find bifurcation points. In this figure, $r_{A}$ is a bifurcation point, and there is another bifurcation point on the interval $(r_{B},r_{C})$. However, on $r_{D}$ we cannot ensure the existence of a bifurcation point, as $h$ does not cross the eigenvalue $\hat{\mu}_{i+2}$.}
\end{figure}

\begin{remark}
As we can see in previous proof, a point of bifurcation always appears when the function $h(r)=(n-1)\left(\dot{\alpha}^2(r)-\ddot{\alpha}(r)\alpha(r)\right)$ cross the value of an eigenvalue of the Laplacian $\Delta$ (see Figure \ref{fig1}), obtaining a weaker condition than (\ref{cond1}) to ensure the existence of bifurcation points. 
\end{remark}

\begin{corolario}
If $r_2=\infty$ and $\alpha(r)=Cr^k + O(r^{k-1})$ with $k>1$ and $C\in\mathbb{R}\slash\left\{0\right\}$, then there are infinitely many bifurcation points.
\end{corolario}
\Demo This is just a direct computation showing that \[\lim_{r\rightarrow \infty}\left( \dot{\alpha}^2(r)-\ddot{\alpha}(r)\alpha(r)\right)=\lim_{r\rightarrow \infty}\left(kC^2r^{2(k-1)}+O(r^{2k-3})\right)=\infty.\]
\cvd\\

\smallskip

Finally, we devote the rest of this section to give a more geometrical viewpoint for previous results. For this, let us recall the formulae for the mean and scalar curvatures on warped products. The mean curvature of $\bord{2}$ is given by

\begin{align}\label{Hwarped}
H(r)=-(n-1)\frac{\dot{\alpha}(r)}{\alpha(r)}
\end{align}
 while the scalar curvature in a point $(r,x)\in (r_1,r_2)\times P$ is given by

\begin{align*}
R{(r,x)}=\frac{1}{\alpha^2(r)}\left(R(x)-(n-1)\left((n-2)\dot{\alpha}^2(r) + 2\ddot{\alpha}(r)\alpha(r)\right)\right)
\end{align*}
where $R(x)$ denote the scalar curvature of the manifold $P$.
Taking into account the expression of $H(r)$, we are able to rewrite previous expression in the following way:

\begin{align*}
\begin{array}{rl} R(r,x)= & \dfrac{R(x)}{\alpha^2(r)}+2(n-1)\dfrac{\dot{\alpha}^2(r)-\ddot{\alpha}(r)\alpha(r)}{\alpha^2(r)}-n(n-1)\dfrac{\dot{\alpha}^2(r)}{\alpha^2(r)}\\\\ =&
\dfrac{R(x)}{\alpha^2(r)}+2(n-1)\dfrac{\dot{\alpha}^2(r)-\ddot{\alpha}(r)\alpha(r)}{\alpha^2(r)} - \dfrac{n}{n-1}H(r)^2.
\end{array}
\end{align*}

In particular,

\[
\frac{\dot{\alpha}^2(r)-\ddot{\alpha}(r)\alpha(r)}{\alpha^2(r)}=\frac{1}{2(n-1)}\left(R(r,x)-\frac{R(x)}{\alpha^2(r)}+\frac{n}{n-1}H(r)^2 \right)
\]
and so, substituting such expression in (\ref{eqbif}), we obtain the following 

\[
{\mathcal A}_{r}(f)=\Delta f -\frac{1}{2}\left(R(r,x)-\frac{R(x)}{\alpha^2(r)}+\frac{n}{n-1}H(r)^2\right) f.
\]

With this expression, we are also able to rewrite the eigenvalues of the operator ${\mathcal A}_{r}$ (compare with \eqref{valfun} and \eqref{valfun1})

\[
\overline{\mu_{i}}(r)=\frac{1}{\alpha^2(r)}\left(\hat{\mu}_{i}+\frac{R(x)}{2}-\frac{\alpha^2(r)}{2}\left(R(r,x) + \frac{n}{n-1} H(r)^2 \right) \right), \;\textrm{ with }i\in\mathbb{N}.
\]

Using these expressions, we can translate Propositions \ref{prop2.3} and \ref{prop5.1} in more geometrical terms. In particular, we deduce the following two corollaries.

\begin{corolario}\label{corsc}
Assume that $R(r,x), R(x) $ and $H(r)$ satisfy
\[
\alpha^2(r)\left(R(r,x) + \frac{n}{n-1}H(r)^2\right)<2\hat{\mu}_1+R(x)
\]
for all $r$ and all $x\in P$. Then, there is no bifurcation point associated to our variational problem.

\end{corolario}

\begin{example}
We consider the following warped product \[\left((0,r_2)\times \mathbb{S}^{n-1}, \,\!dr^2+\sinh^2r\,g^{\mathbb{S}^{n-1}}\!\right).\]
Fixing the set $\Sigma_1=\left\{0\right\}\times\mathbb{S}^{n-1}$ it is possible to visualize the CMC hypersurfaces $\left\{r\right\}\times\mathbb{S}^{n-1}$ as the other set $\Sigma_2$ forming the boundary of the elements of an one-parameter family $\left\{\opens_r\right\}_r$ in this warped product. By previous corollary and the fact that $R(r,x)=-n(n-1)$, $R(x)=(n-1)(n-2)$, $H(r)^2=(n-1)^2\coth^2 r$ and $\hat{\mu}_1=n-1$ we conclude that this family is rigid.
\end{example}
\begin{example}
Another interesting example occurs in the spatial fiber of the Sitter cusp spacetime (see \cite{H} for more details about this spacetime), that is, the following warped product \[\left(
(r_1,r_2)\times \mathbb{T}^{n-1},\,dr^2+e^{2r}\,g^{\mathbb{T}^{n-1}}\!\right),\]
where $ \mathbb{T}^{n-1}$ is the $(n-1)$-dimensional flat torus.
Note that if we fix a set $\Sigma_1=\left\{r_1\right\}\times\mathbb{T}^{n-1}$ for some $r_1\in\mathbb{R}$ and if we consider the CMC hypersurfaces $\left\{r\right\}\times\mathbb{T}^{n-1}$ as the set $\bord{2}$ constituting together with $\Sigma_1$ the boundary of the elements of an one-parameter family $\left\{\opens_r\right\}_r$  we have that $R(r,x)\!=-n(n-1)$, $R(x)=0$, $H(r)^2=(n-1)^2$ and as $\mathbb{T}^{n-1}$ is a compact manifold is clear that $\hat{\mu}_1>0$. Moreover, from \cite[Theorem B.I.2]{BGM} we have that its spectrum is given by $\left\{\hat{\mu}=4\pi^2\left|y\right|^2\right\}$ where $y$ is the closest element to the unity in $\Gamma^*$ (the dual lattice of $\Gamma$ that define the torus in $\mathbb{R}^{n-1}$). Thus, using Corollary \ref{corsc} we conclude that this family has no bifurcation point.
\end{example}

\begin{corolario}
If
\[
\lim_{r\rightarrow r_2}\alpha^2(r)\left(R(r,x) +\frac{n}{n-1}H(r)^2\right)=\infty,
\]
then there exist infinitely many bifurcation points for the family $\left\{\Omega_{r}\right\}_r$ associated to the functional $\biffun{\lambda(r)}$ and the warped product space $(r_1,r_2)\times P$ where $P$ is a $(n-1)$-dimensional closed manifold.
\end{corolario}

\smallskip

\section{\textbf{An application on the de Sitter and Anti-de Sitter Schwarzschild models.} }\label{App}

In this section, we will apply previous results to the following Schwarzschild models: consider a $3$-dimensional Riemannian manifold $M= I\times \mathbb{S}^2$ endowed with the metric

\begin{align}\label{metricschwgen}
g_{K,E}=\psi_{K,E}(r)^{-2}dr^2+r^2\left(d\theta^2+sin^2\theta d\phi^2 \right),
\end{align}
where $K\neq0$ and $\E$ are constants, $\psi_{K,E}(r)=\sqrt{1-2K/r+\E r^2}$ and
$I$ is a maximal connected and open interval where $\psi_{K,E}$ is well defined. For instance, for $K>0$ and $E\geq 0$, $I=(\widehat{r},\infty)$ (with $\widehat{r}>0$) while for $K>0$ and $E<0$ the interval $I=(\widehat{r},0)$ where, in both cases, $\widehat{r}$ is a zero of the function $\psi_{K,E}$.

The two parameters $K$ and $\E$ represent \textit{the black hole mass} (and so, it is usually considered positive) and \textit{the cosmological constant} respectively. When $E=0$, we obtain the classical model of the spatial fiber of the \textit{Schwarzschild spacetime}. For $E$ negative (resp. positive) we have the Riemannian model of the so-called de \textit{Sitter} (resp. \textit{Anti-de Sitter}) \textit{Schwarzschild spacetime}.

\smallskip

It follows directly that the variable change
\begin{equation}\label{varchan}ds=-\psi_{K,E}(r)^{-1}dr
\end{equation}
transform the metric $g_{K,E}$ into the warped metric

\begin{equation}\label{varchanmet}
g_{K,E}=ds^2+r(s)^2\left(d\theta^2+sin^2\theta \,d\phi^2 \right)
\end{equation}

\noindent making the results of previous section applicable. In this case, $\alpha(s)=r(s)$ and the CMC surfaces $\left\{s\right\}\times \mathbb{S}^2$ are $2$-spheres of radius $r(s)$, which can be seen as the boundaries of open subsets $\Omega_s$ that belong to a one-parameter family of $(M^{3},g_{K,E})$. Moreover, the eigenvalues for the $2$-spheres are well known (recall \cite[Proposition C.I.1]{BGM}), and satisfy that $\hat{\mu}_i=i(i+1)$ with $i\in \mathbb{N}_0$.

\smallskip

Let us consider two relevant cases (the rest will follow analogously). Consider first that both parameters $K$ and $E$ are positive (i.e., we consider a Anti-de Sitter Schwarzschild model) and take the interval $I=(\widehat{r},\infty)$ where $\widehat{r}$ is the only positive zero of $\psi_{K,E}$. Although the expression of $g_{K,E}$ becomes singular on $\widehat{r}$, it can be proved that the metric $g_{K,E}$ extends to a smooth Riemannian metric on $M=[\widehat{r},\infty)\times \mathbb{S}^2$. Finally, by the variable change on \eqref{varchan}, $M=[\widehat{s},\infty)\times \mathbb{S}^2$ and the metric $g_{K,E}$ takes the form on \eqref{varchanmet}. Then, we can state the following result

\begin{corolario}
Consider the spatial fiber of Anti-de Sitter Schwarzschild spacetime $(M^3,g_{K,E})$ (with $K,E>0$) as a warped product space. The family $\left\{\Omega_s\right\}_{s\in (\widehat{s},\infty)}$ where $\Omega_s=(\widehat{s},s)\times \mathbb{S}^2$ is locally rigid on all $s>\widehat{s}$.
\end{corolario}
\Demo
In order to apply Proposition \ref{prop2.3}, let us estimate the following:

\begin{align}\label{rl}
\dot{r}^2(s)-\ddot{r}(s)r(s) =1-\dfrac{3K}{r(s)}
\end{align}

where we have used that $\dot{r}(s)=-\sqrt{1-\dfrac{2K}{r(s)}+\E r(s)^2}$.\\\\
Now, bearing in mind that $\hat{\mu}_{1}=2$ and $n=3$, we obtain that

\[
\dot{r}^2(s)-\ddot{r}(s)r(s)=1-\frac{3K}{r(s)}
<\frac{\hat{\mu}_1}{n-1}=1
\]
for $r(s)>\widehat{r}>0$ and $K>0$. So, Proposition \ref{prop2.3} applies showing that there are no bifurcation points. \cvd
\\
\smallskip

Another way to argue the rigidity of this family $\left\{\Omega_s\right\}_s$ is to use Corollary \ref{corsc} and the fact $R(r,x)=-6\E$, $H(r)^2=4/r^2\psi_{K,E}(r)^2=4/r^2-8K/r^3+4\E$ and $R(x)=2$.\\

Now, let us consider the case when the cosmological constant $E$ is negative, and let us take $I=[\widehat{r},0)$ (again, the Riemannian metric extends smoothly). Take the variable change \eqref{varchan} in such a way that $M=[\widehat{s},0)\times \mathbb{S}^2$. In this case, the family $\{\Omega_{s}\}_{s\in(\widehat{s},0)}$ is no longer locally rigid, as we can see on the following result

\begin{corolario}
Consider the spatial fiber of de Sitter Schwarzschild spacetime $(M^3,g_{K,E})$ (i.e., with $E<0$ and $K>0$) as a warped product space. Then the family $\left\{\Omega_s\right\}_{s\in (\widehat{s},s)}$ where $\Omega_s=(\widehat{s},s)\times \mathbb{S}^2$ has infinitely many bifurcation points.
\end{corolario}
\Demo
Observe that
\begin{equation}
\lim_{s\rightarrow 0^-}\dot{r}^2(s)-\ddot{r}(s)r(s)=\lim_{r\rightarrow 0^-}1-\frac{3K}{r}=\infty.
\end{equation}
 Then, by Proposition \ref{prop5.1} we conclude that the family $\left\{\Omega_s\right\}_s$ has infinitely many bifurcation points.
\cvd

\bibliographystyle{amsplain}

\end{document}